\documentclass[10pt,draft]{article}
\usepackage{amsfonts,amsmath,amssymb,cite}
\numberwithin{equation}{section}
\begin{document}
\title{On parameter derivatives of the associated Legendre function
of the first kind (with applications to the construction of the
associated Legendre function of the second kind of integer degree and
order)}
\author{Rados{\l}aw Szmytkowski \\*[3ex]
Atomic Physics Division,
Department of Atomic Physics and Luminescence, \\
Faculty of Applied Physics and Mathematics,
Gda{\'n}sk University of Technology, \\
Narutowicza 11/12, PL 80--233 Gda{\'n}sk, Poland \\
email: radek@mif.pg.gda.pl}
\date{\today}
\maketitle
\begin{abstract}
A relationship between partial derivatives of the associated Legendre
function of the first kind with respect to its degree, $[\partial
P_{\nu}^{m}(z)/\partial\nu]_{\nu=n}$, and to its order, $[\partial
P_{n}^{\mu}(z)/\partial\mu]_{\mu=m}$, is established for
$m,n\in\mathbb{N}$. This relationship is used to deduce four new
closed-form representations of $[\partial
P_{\nu}^{m}(z)/\partial\nu]_{\nu=n}$ from those found recently for
$[\partial P_{n}^{\mu}(z)/\partial\mu]_{\mu=m}$ by the present author
[R.~Szmytkowski, J.\ Math.\ Chem.\ 46 (2009) 231]. Several new
expressions for the associated Legendre function of the second kind
of integer degree and order, $Q_{n}^{m}(z)$, suitable for numerical
purposes, are also derived.
\end{abstract}
\maketitle
%
%
\section{Introduction}
\label{I}
\setcounter{equation}{0}
Recently, an interest has arisen in the derivation of closed-form
expressions for parameter derivatives of the associated Legendre
functions \cite{Bryc09,Cohl09,Szmy09a,Szmy09b}. In Refs.\
\cite{Szmy09a,Szmy09b}, the present author has extensively studied
such derivatives for the associated Legendre function function of the
first kind, $P_{\nu}^{\mu}(z)$, in the case when one of the
parameters is a fixed integer. In particular, in Ref.\
\cite{Szmy09a}, using finite-sum expressions for $P_{n}^{\mu}(z)$,
with $n\in\mathbb{N}$, we have found the following two
representations of $[\partial P_{n}^{\mu}(z)/\partial\mu]_{\mu=m}$,
with $m\in\mathbb{N}$:
\begin{eqnarray}
\frac{\partial P_{n}^{\mu}(z)}{\partial\mu}\bigg|_{\mu=m} 
&=& \frac{1}{2}P_{n}^{m}(z)\ln\frac{z+1}{z-1}
\nonumber \\
&& +\,(-)^{m}\left(\frac{z+1}{z-1}\right)^{m/2}
\sum_{k=0}^{m-1}(-)^{k}\frac{(k+n)!(m-k-1)!}{k!(n-k)!}
\left(\frac{z-1}{2}\right)^{k}
\nonumber \\
&& +\,\left(\frac{z^{2}-1}{4}\right)^{m/2}
\sum_{k=0}^{n-m}\frac{(k+n+m)!\psi(k+1)}{k!(k+m)!(n-m-k)!}
\left(\frac{z-1}{2}\right)^{k}
\qquad (0\leqslant m\leqslant n)
\nonumber \\
&&
\label{1.1}
\end{eqnarray}
and
\begin{eqnarray}
\frac{\partial P_{n}^{\mu}(z)}{\partial\mu}\bigg|_{\mu=m} 
&=& \frac{1}{2}P_{n}^{m}(z)\ln\frac{z+1}{z-1}
+[\psi(n+m+1)+\psi(n-m+1)]P_{n}^{m}(z)
\nonumber \\
&& -\,(-)^{n}\frac{(n+m)!}{(n-m)!}\left(\frac{z+1}{z-1}\right)^{m/2}
\sum_{k=0}^{n}(-)^{k}
\frac{(k+n)!\psi(k+m+1)}{k!(k+m)!(n-k)!}
\left(\frac{z+1}{2}\right)^{k}
\nonumber \\
&& (0\leqslant m\leqslant n).
\label{1.2}
\end{eqnarray}
In turn, in Ref.\ \cite{Szmy09b}, using contour-integral
representations of $\partial P_{\nu}^{m}(z)/\partial\nu$, we have
arrived, among others, at the following three expressions for
$[\partial P_{\nu}^{m}(z)/\partial\nu]_{\nu=n}$, again with
$m,n\in\mathbb{N}$:
\begin{eqnarray}
\frac{\partial P_{\nu}^{m}(z)}{\partial\nu}\bigg|_{\nu=n}
&=& P_{n}^{m}(z)\ln\frac{z+1}{2}
-[\psi(n+1)+\psi(n-m+1)]P_{n}^{m}(z)
\nonumber \\
&& +\,\left(\frac{z^{2}-1}{4}\right)^{m/2}
\sum_{k=0}^{n-m}\frac{(k+n+m)!\psi(k+n+m+1)}{k!(k+m)!(n-m-k)!}
\left(\frac{z-1}{2}\right)^{k}
\nonumber \\
&& +\,\frac{(n+m)!}{(n-m)!}\left(\frac{z-1}{z+1}\right)^{m/2}
\sum_{k=0}^{n}\frac{(k+n)!\psi(k+n+1)}{k!(k+m)!(n-k)!}
\left(\frac{z-1}{2}\right)^{k}
\nonumber \\
&& (0\leqslant m\leqslant n),
\label{1.3}
\end{eqnarray}
\begin{eqnarray}
\frac{\partial P_{\nu}^{m}(z)}{\partial\nu}\bigg|_{\nu=n}
&=& P_{n}^{m}(z)\ln\frac{z+1}{2}
+[\psi(n+1)-\psi(n-m+1)]P_{n}^{m}(z)
\nonumber \\
&& -\,(-)^{n}\frac{(n+m)!}{(n-m)!}
\left(\frac{z^{2}-1}{4}\right)^{-m/2}
\sum_{k=0}^{m-1}\frac{(k+n-m)!(m-k-1)!}{k!(n+m-k)!}
\left(\frac{z+1}{2}\right)^{k}
\nonumber \\
&& +\,(-)^{n+m}\left(\frac{z^{2}-1}{4}\right)^{m/2}
\sum_{k=0}^{n-m}(-)^{k}\frac{(k+n+m)!}{k!(k+m)!(n-m-k)!}
\nonumber \\
&& \quad \times[\psi(k+n+m+1)-\psi(k+m+1)]
\left(\frac{z+1}{2}\right)^{k}
\nonumber \\
&& +\,(-)^{n}\frac{(n+m)!}{(n-m)!}\left(\frac{z+1}{z-1}\right)^{m/2}
\sum_{k=0}^{n}(-)^{k}\frac{(k+n)!}{k!(k+m)!(n-k)!}
\nonumber \\
&& \quad \times[\psi(k+n+1)-\psi(k+1)]\left(\frac{z+1}{2}\right)^{k}
\qquad (0\leqslant m\leqslant n)
\label{1.4}
\end{eqnarray}
and
\begin{eqnarray}
\frac{\partial P_{\nu}^{m}(z)}{\partial\nu}\bigg|_{\nu=n} 
&=& P_{n}^{m}(z)\ln\frac{z+1}{2}
+[\psi(n+m+1)-\psi(n+1)]P_{n}^{m}(z)
\nonumber \\
&& -\,(-)^{n+m}\left(\frac{z-1}{z+1}\right)^{m/2}\sum_{k=0}^{m-1}
\frac{(k+n)!(m-k-1)!}{k!(n-k)!}\left(\frac{z+1}{2}\right)^{k}
\nonumber \\
&& +\,(-)^{n+m}\left(\frac{z^{2}-1}{4}\right)^{m/2}
\sum_{k=0}^{n-m}(-)^{k}\frac{(k+n+m)!}{k!(k+m)!(n-m-k)!}
\nonumber \\
&& \quad\times[\psi(k+n+m+1)-\psi(k+1)]
\left(\frac{z+1}{2}\right)^{k}
\nonumber \\
&& +\,(-)^{n}\frac{(n+m)!}{(n-m)!}\left(\frac{z+1}{z-1}\right)^{m/2}
\sum_{k=0}^{n}(-)^{k}\frac{(k+n)!}{k!(k+m)!(n-k)!}
\nonumber \\
&& \quad\times[\psi(k+n+1)-\psi(k+m+1)]
\left(\frac{z+1}{2}\right)^{k}
\qquad (0\leqslant m\leqslant n).
\label{1.5}
\end{eqnarray}
In the above equations, and in what follows, $\psi(\zeta)$ is the
digamma function defined as
\begin{equation}
\psi(\zeta)=\frac{1}{\Gamma(\zeta)}
\frac{\mathrm{d}\Gamma(\zeta)}{\mathrm{d}\zeta}.
\label{1.6}
\end{equation}

In the present paper, we shall pursue further the subject of
derivation of closed-form expressions for parameter derivatives of
the associated Legendre function of the first kind. First, in Section
\ref{II}, we shall show that there exists a simple relationship
between the derivatives $[\partial
P_{\nu}^{m}(z)/\partial\nu]_{\nu=n}$ and $[\partial
P_{n}^{\mu}(z)/\partial\mu]_{\mu=m}$, both with $m,n\in\mathbb{N}$.
Next, in Section \ref{III}, we shall use this relationship, in
conjunction with the formulas (\ref{1.1}) and (\ref{1.2}), to derive
four further representations of $[\partial
P_{\nu}^{m}(z)/\partial\nu]_{\nu=n}$, two of them involving sums of
powers of $(z+1)/2$ and the remaining two --- sums of powers of
$(z-1)/2$. Interestingly, each out of these four representations
contains only \emph{two\/} sums. Therefore, the two new expressions
for $[\partial P_{\nu}^{m}(z)/\partial\nu]_{\nu=n}$ containing sums
of powers of $(z+1)/2$ appear to be markedly simpler than the
previously derived representations (\ref{1.4}) and (\ref{1.5}), while
the two new expressions involving sums of powers of $(z-1)/2$ are of
the same degree of complexity as the representation (\ref{1.3}). In
addition, as a by-product, we shall obtain in that section two useful
finite-sum representations of $[\psi(n+m+1)-\psi(n+1)]P_{n}^{m}(z)$.
In the final Section \ref{IV}, we shall exploit the results of
Section \ref{III} to find some new representations of the associated
Legendre function of the second kind of integer degree and order,
$Q_{n}^{m}(z)$, suitable for use for numerical purposes in various
parts of the complex $z$-plane.

Throughout the paper, we shall be adopting the standard convention
according to which $z\in\mathbb{C}\setminus[-1,1]$, with the phases
restricted by
\begin{equation}
-\pi<\arg(z)<\pi, \qquad -\pi<\arg(z\pm1)<\pi
\label{1.7}
\end{equation}
(this corresponds to drawing a cut in the $z$-plane along the real
axis from $z=-\infty$ to $z=+1$), hence,
\begin{equation}
-z=\mathrm{e}^{\mp\mathrm{i}\pi}z
\qquad
-z+1=\mathrm{e}^{\mp\mathrm{i}\pi}(z-1)
\qquad
-z-1=\mathrm{e}^{\mp\mathrm{i}\pi}(z+1)
\qquad (\arg(z)\gtrless0).
\label{1.8}
\end{equation}
Also, it will be implicit that $x\in[-1,1]$, $\mu,\nu\in\mathbb{C}$
and $k,m,n\in\mathbb{N}$. Finally, it will be understood that if the
upper limit of a sum is less by unity than the lower one, then the
sum vanishes identically. The associated Legendre functions of the
first and the second kinds used in the paper are those of Hobson
\cite{Hobs31} (cf.\ also Refs.\ \cite{Robi57,Robi58,Robi59}).

Before proceeding to the matter, we emphasize that the results
obtained in the present paper are interesting not only for their own
(mathematical) sake. The derivatives $\partial
P_{\nu}^{m}(z)/\partial\nu$ and $\partial P_{\nu}^{m}(x)/\partial\nu$
are met in solutions of some boundary value problems of theoretical
acoustics, electromagnetism, heat conduction and other branches of
theoretical physics and applied mathematics (for a list of references
illustrating this statement, see Ref.\ \cite{Szmy09b}). Physical
applications of the associated Legendre functions of the second kind
of integer degree and order, $Q_{n}^{m}(z)$ and $Q_{n}^{m}(x)$, (for
which, recall, several new explicit expressions are derived in
Section \ref{IV}) are even more abundant \cite{Mors53,Moon61}.
%
%
\section{A relationship between $[\partial
P_{\nu}^{m}(z)/\partial\nu]_{\nu=n}$ and $[\partial
P_{n}^{\mu}(z)/\partial\mu]_{\mu=m}$}
\label{II}
\setcounter{equation}{0}
The departure point for our considerations in this section is the
following Rodrigues-type formula, due to Barnes \cite{Barn07}, for
the associated Legendre function of the first kind when the sum of
its degree and its order is a non-negative integer:
\begin{equation}
P_{\nu}^{m-\nu}(z)=\frac{1}{2^{\nu}\Gamma(\nu+1)}
(z^{2}-1)^{(m-\nu)/2}
\frac{\mathrm{d}^{m}(z^{2}-1)^{\nu}}{\mathrm{d}z^{m}}.
\label{2.1}
\end{equation}
In terms of the Jacobi polynomial
\begin{equation}
P_{n}^{(\alpha,\beta)}(z)=\frac{1}{2^{n}n!}
(z-1)^{-\alpha}(z+1)^{-\beta}
\frac{\mathrm{d}^{n}}{\mathrm{d}z^{n}}
\left[(z-1)^{n+\alpha}(z+1)^{n+\beta}\right]
\qquad (\alpha,\beta\in\mathbb{C}),
\label{2.2}
\end{equation}
Eq.\ (\ref{2.1}) may be rewritten as
\begin{equation}
P_{\nu}^{m-\nu}(z)=\frac{m!}{\Gamma(\nu+1)}
\left(\frac{z^{2}-1}{4}\right)^{(\nu-m)/2}P_{m}^{(\nu-m,\nu-m)}(z).
\label{2.3}
\end{equation}
(With no doubt, the reader has immediately realized that the Jacobi
polynomial appearing on the right-hand side of Eq.\ (\ref{2.3}) is a
multiple of the Gegenbauer polynomial $C_{m}^{(\nu-m+1/2)}(z)$.
However, we shall not make any use of this fact here.)
Differentiation of Eq.\ (\ref{2.3}) with respect to $\nu$, followed
by setting $\nu=n$, yields
\begin{eqnarray}
\frac{\partial P_{\nu}^{m-n}(z)}{\partial\nu}\bigg|_{\nu=n}
-\frac{\partial P_{n}^{\mu}(z)}{\partial\mu}\bigg|_{\mu=m-n}
&=& \frac{1}{2}P_{n}^{m-n}(z)\ln\frac{z^{2}-1}{4}
-\psi(n+1)P_{n}^{m-n}(z)
\nonumber \\
&& +\frac{m!}{n!}\left(\frac{z^{2}-1}{4}\right)^{(n-m)/2}
\frac{\partial P_{m}^{(\lambda,\lambda)}(z)}
{\partial\lambda}\bigg|_{\lambda=n-m}.
\label{2.4}
\end{eqnarray}
The replacement of $m$ by $n+m$ results in the relationship
\begin{eqnarray}
\frac{\partial P_{\nu}^{m}(z)}{\partial\nu}\bigg|_{\nu=n}
-\frac{\partial P_{n}^{\mu}(z)}{\partial\mu}\bigg|_{\mu=m}
&=& \frac{1}{2}P_{n}^{m}(z)\ln\frac{z^{2}-1}{4}
-\psi(n+1)P_{n}^{m}(z)
\nonumber \\
&& +\frac{(n+m)!}{n!}\left(\frac{z^{2}-1}{4}\right)^{-m/2}
\frac{\partial P_{n+m}^{(\lambda,\lambda)}(z)}
{\partial\lambda}\bigg|_{\lambda=-m}.
\label{2.5}
\end{eqnarray}
If in Eq.\ (\ref{2.5}) one exploits the following two explicit
representations of the Jacobi polynomial $P_{n}^{(\alpha,\beta)}(z)$:
\begin{equation}
P_{n}^{(\alpha,\beta)}(z)
=\frac{\Gamma(n+\alpha+1)}{\Gamma(n+\alpha+\beta+1)}
\sum_{k=0}^{n}\frac{\Gamma(k+n+\alpha+\beta+1)}
{k!(n-k)!\Gamma(k+\alpha+1)}\left(\frac{z-1}{2}\right)^{k},
\label{2.6}
\end{equation}
\begin{equation}
P_{n}^{(\alpha,\beta)}(z)
=(-)^{n}\frac{\Gamma(n+\beta+1)}{\Gamma(n+\alpha+\beta+1)}
\sum_{k=0}^{n}(-)^{k}\frac{\Gamma(k+n+\alpha+\beta+1)}
{k!(n-k)!\Gamma(k+\beta+1)}\left(\frac{z+1}{2}\right)^{k},
\label{2.7}
\end{equation}
after making use of Eq.\ (\ref{2.3}) and of the easily provable
relation 
\begin{equation}
\lim_{\lambda\to-m}\frac{\psi(k+\lambda+1)}{\Gamma(k+\lambda+1)}
=(-)^{k+m}(m-k-1)!
\qquad (0\leqslant k\leqslant m-1),
\label{2.8}
\end{equation}
one obtains
\begin{eqnarray}
&& \frac{\partial P_{\nu}^{m}(z)}{\partial\nu}\bigg|_{\nu=n}
-\frac{\partial P_{n}^{\mu}(z)}{\partial\mu}\bigg|_{\mu=m}
= \frac{1}{2}P_{n}^{m}(z)\ln\frac{z^{2}-1}{4}
-2\psi(n-m+1)P_{n}^{m}(z)
\nonumber \\
&& \hspace*{5em}
-\,(-)^{m}\frac{(n+m)!}{(n-m)!}
\left(\frac{z^{2}-1}{4}\right)^{-m/2}
\sum_{k=0}^{m-1}(-)^{k}\frac{(k+n-m)!(m-k-1)!}{k!(n+m-k)!}
\left(\frac{z-1}{2}\right)^{k}
\nonumber \\
&& \hspace*{5em}
+\,\frac{(n+m)!}{(n-m)!}\left(\frac{z-1}{z+1}\right)^{m/2}
\sum_{k=0}^{n}\frac{(k+n)!}{k!(k+m)!(n-k)!}
\nonumber \\
&& \hspace*{6em}
\times[2\psi(k+n+1)-\psi(k+1)]\left(\frac{z-1}{2}\right)^{k}
\qquad (0\leqslant m\leqslant n)
\label{2.9}
\end{eqnarray}
and
\begin{eqnarray}
&& \frac{\partial P_{\nu}^{m}(z)}{\partial\nu}\bigg|_{\nu=n}
-\frac{\partial P_{n}^{\mu}(z)}{\partial\mu}\bigg|_{\mu=m}
=\frac{1}{2}P_{n}^{m}(z)\ln\frac{z^{2}-1}{4}
-2\psi(n-m+1)P_{n}^{m}(z)
\nonumber \\
&& \hspace*{5em}
-\,(-)^{n}\frac{(n+m)!}{(n-m)!}\left(\frac{z^{2}-1}{4}\right)^{-m/2}
\sum_{k=0}^{m-1}\frac{(k+n-m)!(m-k-1)!}{k!(n+m-k)!}
\left(\frac{z+1}{2}\right)^{k}
\nonumber \\
&& \hspace*{5em}
+\,(-)^{n}\frac{(n+m)!}{(n-m)!}\left(\frac{z+1}{z-1}\right)^{m/2}
\sum_{k=0}^{n}(-)^{k}\frac{(k+n)!}{k!(k+m)!(n-k)!}
\nonumber \\
&& \hspace*{6em}
\times[2\psi(k+n+1)-\psi(k+1)]\left(\frac{z+1}{2}\right)^{k}
\qquad (0\leqslant m\leqslant n), 
\label{2.10}
\end{eqnarray}
respectively.
%
%
\section{Some new representations of $[\partial 
P_{\nu}^{m}(z)/\partial\nu]_{\nu=n}$ with \mbox{$0\leqslant
m\leqslant n$}}
\label{III} 
\setcounter{equation}{0}
In this section, we shall use the results of Section \ref{II} to
provide several expressions for the derivative $[\partial
P_{\nu}^{m}(z)/\partial\nu]_{\nu=n}$, which differ from these given
in Eqs.\ (\ref{1.3})--(\ref{1.5}).

The first from among these expressions for $[\partial
P_{\nu}^{m}(z)/\partial\nu]_{\nu=n}$ follows if one plugs the
representation (\ref{1.2}) of $[\partial
P_{n}^{\mu}(z)/\partial\mu]_{\mu=m}$ into the relationship
(\ref{2.10}). The result
\begin{eqnarray}
\frac{\partial P_{\nu}^{m}(z)}{\partial\nu}\bigg|_{\nu=n} 
&=& P_{n}^{m}(z)\ln\frac{z+1}{2}
+[\psi(n+m+1)-\psi(n-m+1)]P_{n}^{m}(z)
\nonumber \\
&& -\,(-)^{n}\frac{(n+m)!}{(n-m)!}
\left(\frac{z^{2}-1}{4}\right)^{-m/2}\sum_{k=0}^{m-1}
\frac{(k+n-m)!(m-k-1)!}{k!(n+m-k)!}\left(\frac{z+1}{2}\right)^{k}
\nonumber \\
&& +\,(-)^{n}\frac{(n+m)!}{(n-m)!}\left(\frac{z+1}{z-1}\right)^{m/2}
\sum_{k=0}^{n}(-)^{k}\frac{(k+n)!}{k!(k+m)!(n-k)!}
\nonumber \\
&& \quad\times[2\psi(k+n+1)-\psi(k+m+1)-\psi(k+1)]
\left(\frac{z+1}{2}\right)^{k}
\qquad (0\leqslant m\leqslant n)
\nonumber \\
&&
\label{3.1}
\end{eqnarray}
is seen to be much simpler than either of the representations
(\ref{1.4}) or (\ref{1.5}). An infinite variety of other
representations of $[\partial P_{\nu}^{m}(z)/\partial\nu]_{\nu=n}$,
involving sums of powers of $(z+1)/2$, may be obtained by taking
linear combinations, with coefficients such that their sum is unity,
of the expressions in Eqs.\ (\ref{1.4}), (\ref{1.5}) and (\ref{3.1}).
For instance, multiplying Eq.\ (\ref{3.1}) by $-1$ and adding to the
sum of Eqs.\ (\ref{1.4}) and (\ref{1.5}) leads to another remarkably
simple formula
\begin{eqnarray}
\frac{\partial P_{\nu}^{m}(z)}{\partial\nu}\bigg|_{\nu=n} 
&=& P_{n}^{m}(z)\ln\frac{z+1}{2}
\nonumber \\
&& -\,(-)^{n+m}\left(\frac{z-1}{z+1}\right)^{m/2}\sum_{k=0}^{m-1}
\frac{(k+n)!(m-k-1)!}{k!(n-k)!}\left(\frac{z+1}{2}\right)^{k}
\nonumber \\
&& +\,(-)^{n+m}\left(\frac{z^{2}-1}{4}\right)^{m/2}
\sum_{k=0}^{n-m}(-)^{k}\frac{(k+n+m)!}{k!(k+m)!(n-m-k)!}
\nonumber \\
&& \quad\times[2\psi(k+n+m+1)-\psi(k+m+1)-\psi(k+1)]
\left(\frac{z+1}{2}\right)^{k}
\nonumber \\
&& (0\leqslant m\leqslant n).
\label{3.2}
\end{eqnarray}
It is seen that for $m=0$ both the representations (\ref{3.1}) and
(\ref{3.2}) of $[\partial P_{\nu}^{m}(z)/\partial\nu]_{\nu=n}$ reduce
to the formula
\begin{eqnarray}
\frac{\partial P_{\nu}(z)}{\partial\nu}\bigg|_{\nu=n}
&=& P_{n}(z)\ln\frac{z+1}{2}
+2\sum_{k=0}^{n}(-)^{k+n}\frac{(k+n)!}{(k!)^{2}(n-k)!}
[\psi(k+n+1)-\psi(k+1)]\left(\frac{z+1}{2}\right)^{k},
\nonumber \\
&&
\label{3.3}
\end{eqnarray}
found by the present author in Ref.\ \cite[Section 5.2.7]{Szmy06}
(cf.\ also Ref.\ \cite{Szmy07}).

From the above findings, one may deduce two interesting and, as we
shall see in a moment, useful identities involving the function
$P_{n}^{m}(z)$. If we equate the right-hand sides of Eqs.\
(\ref{1.4}) and (\ref{3.1}), this results in the first of these
relations:
\begin{eqnarray}
[\psi(n+m+1)-\psi(n+1)]P_{n}^{m}(z)
&=& (-)^{n+m}\left(\frac{z^{2}-1}{4}\right)^{m/2}
\sum_{k=0}^{n-m}(-)^{k}\frac{(k+n+m)!}{k!(k+m)!(n-m-k)!}
\nonumber \\
&& \quad
\times[\psi(k+n+m+1)-\psi(k+m+1)]\left(\frac{z+1}{2}\right)^{k}
\nonumber \\
&& -\,(-)^{n}\frac{(n+m)!}{(n-m)!}\left(\frac{z+1}{z-1}\right)^{m/2}
\sum_{k=0}^{n}(-)^{k}\frac{(k+n)!}{k!(k+m)!(n-k)!}
\nonumber \\
&& \quad
\times[\psi(k+n+1)-\psi(k+m+1)]\left(\frac{z+1}{2}\right)^{k}
\nonumber \\
&& (0\leqslant m\leqslant n).
\label{3.4}
\end{eqnarray}
Replacement of $z$ by $-z$ in the above equation, followed by the use
of the well-known property
\begin{equation}
P_{n}^{m}(-z)=(-)^{n}P_{n}^{m}(z)
\qquad (0\leqslant m\leqslant n)
\label{3.5}
\end{equation}
and also of Eq.\ (\ref{1.8}), leads to the second identity:
\begin{eqnarray}
[\psi(n+m+1)-\psi(n+1)]P_{n}^{m}(z)
&=& \left(\frac{z^{2}-1}{4}\right)^{m/2}
\sum_{k=0}^{n-m}\frac{(k+n+m)!}{k!(k+m)!(n-m-k)!}
\nonumber \\
&& \quad
\times[\psi(k+n+m+1)-\psi(k+m+1)]\left(\frac{z-1}{2}\right)^{k}
\nonumber \\
&& -\,\frac{(n+m)!}{(n-m)!}\left(\frac{z-1}{z+1}\right)^{m/2}
\sum_{k=0}^{n}\frac{(k+n)!}{k!(k+m)!(n-k)!}
\nonumber \\
&& \quad
\times[\psi(k+n+1)-\psi(k+m+1)]\left(\frac{z-1}{2}\right)^{k}
\nonumber \\
&& (0\leqslant m\leqslant n).
\label{3.6}
\end{eqnarray}

Playing with Eq.\ (\ref{1.3}) and with the identity (\ref{3.6}), one
may obtain an infinite variety of representations of $[\partial
P_{\nu}^{m}(z)/\partial\nu]_{\nu=n}$ containing sums of powers of
$(z-1)/2$. Two examples of such representations are
\begin{eqnarray}
\frac{\partial P_{\nu}^{m}(z)}{\partial\nu}\bigg|_{\nu=n}
&=& P_{n}^{m}(z)\ln\frac{z+1}{2}
-[\psi(n+m+1)+\psi(n-m+1)]P_{n}^{m}(z)
\nonumber \\
&& +\,\left(\frac{z^{2}-1}{4}\right)^{m/2}
\sum_{k=0}^{n-m}\frac{(k+n+m)!}{k!(k+m)!(n-m-k)!}
\nonumber \\
&& \quad
\times[2\psi(k+n+m+1)-\psi(k+m+1)]\left(\frac{z-1}{2}\right)^{k}
\nonumber \\
&& +\,\frac{(n+m)!}{(n-m)!}\left(\frac{z-1}{z+1}\right)^{m/2}
\sum_{k=0}^{n}\frac{(k+n)!\psi(k+m+1)}{k!(k+m)!(n-k)!}
\left(\frac{z-1}{2}\right)^{k}
\nonumber \\
&& (0\leqslant m\leqslant n)
\label{3.7}
\end{eqnarray}
and
\begin{eqnarray}
\frac{\partial P_{\nu}^{m}(z)}{\partial\nu}\bigg|_{\nu=n}
&=& P_{n}^{m}(z)\ln\frac{z+1}{2}
+[\psi(n+m+1)-2\psi(n+1)-\psi(n-m+1)]P_{n}^{m}(z)
\nonumber \\
&& +\,\left(\frac{z^{2}-1}{4}\right)^{m/2}
\sum_{k=0}^{n-m}\frac{(k+n+m)!\psi(k+m+1)}{k!(k+m)!(n-m-k)!}
\left(\frac{z-1}{2}\right)^{k}
\nonumber \\
&& +\,\frac{(n+m)!}{(n-m)!}\left(\frac{z-1}{z+1}\right)^{m/2}
\sum_{k=0}^{n}\frac{(k+n)!}{k!(k+m)!(n-k)!}
\nonumber \\
&& \quad
\times[2\psi(k+n+1)-\psi(k+m+1)]\left(\frac{z-1}{2}\right)^{k}
\qquad (0\leqslant m\leqslant n).
\label{3.8}
\end{eqnarray}
For $m=0$, both Eqs.\ (\ref{3.7}) and (\ref{3.8}) reduce to the
Schelkunoff's formula \cite{Sche41} (cf.\ also Ref.\ \cite[Section
5.2.6]{Szmy06})
\begin{eqnarray}
\frac{\partial P_{\nu}(z)}{\partial\nu}\bigg|_{\nu=n}
&=& P_{n}(z)\ln\frac{z+1}{2}-2\psi(n+1)P_{n}(z)
+2\sum_{k=0}^{n}\frac{(k+n)!\psi(k+n+1)}{(k!)^{2}(n-k)!}
\left(\frac{z-1}{2}\right)^{k}.
\nonumber \\
&&
\label{3.9}
\end{eqnarray}

From the representations of $[\partial
P_{\nu}^{m}(z)/\partial\nu]_{\nu=n}$ found above, one may construct
counterpart representations for $[\partial
P_{\nu}^{-m}(z)/\partial\nu]_{\nu=n}$, using the relationship
\cite{Szmy09b}
\begin{eqnarray}
\frac{\partial P_{\nu}^{-m}(z)}{\partial\nu}\bigg|_{\nu=n}
&=& \frac{(n-m)!}{(n+m)!}
\frac{\partial P_{\nu}^{m}(z)}{\partial\nu}\bigg|_{\nu=n}
-[\psi(n+m+1)-\psi(n-m+1)]P_{n}^{-m}(z)
\nonumber \\
&& (0\leqslant m\leqslant n)
\label{3.10}
\end{eqnarray}
and the well-known property
\begin{equation}
P_{n}^{m}(z)=\frac{(n+m)!}{(n-m)!}P_{n}^{-m}(z)
\qquad (0\leqslant m\leqslant n).
\label{3.11}
\end{equation}
Moreover, it does not offer any difficulty to derive counterpart
expressions on the cut $x\in[-1,1]$ by using the formulas
\begin{eqnarray}
\frac{\partial P_{\nu}^{\pm m}(x)}{\partial\nu}\bigg|_{\nu=n}
&=& \mathrm{e}^{\pm\mathrm{i}\pi m/2}
\frac{\partial P_{\nu}^{\pm m}(x+\mathrm{i}0)}
{\partial\nu}\bigg|_{\nu=n}
=\mathrm{e}^{\mp\mathrm{i}\pi m/2}
\frac{\partial P_{\nu}^{\pm m}(x-\mathrm{i}0)}
{\partial\nu}\bigg|_{\nu=n}
\nonumber \\
&=& \frac{1}{2}\left[\mathrm{e}^{\pm\mathrm{i}\pi m/2}
\frac{\partial P_{\nu}^{\pm m}(x+\mathrm{i}0)}
{\partial\nu}\bigg|_{\nu=n}+\mathrm{e}^{\mp\mathrm{i}\pi m/2}
\frac{\partial P_{\nu}^{\pm m}(x-\mathrm{i}0)}
{\partial\nu}\bigg|_{\nu=n}\right]
\nonumber \\
&&
\label{3.12}
\end{eqnarray}
together with
\begin{equation}
x+1\pm\mathrm{i}0=1+x,
\qquad 
x-1\pm\mathrm{i}0=\mathrm{e}^{\pm\mathrm{i}\pi}(1-x).
\label{3.13}
\end{equation}

Concluding this section we note that, in principle, the relations
(\ref{2.9}) and (\ref{2.10}) might be also used in the opposite
direction, i.e., to construct representations for $[\partial
P_{n}^{\mu}(z)/\partial\mu]_{\mu=m}$ from those known for $[\partial
P_{\nu}^{m}(z)/\partial\nu]_{\nu=n}$. As it appears, however, that
all expressions for $[\partial P_{n}^{\mu}(z)/\partial\mu]_{\mu=m}$
obtainable in this way are much more complex (and thus potentially
less useful) than these in Eqs.\ (\ref{1.1}) and (\ref{1.2}), we do
not present them here.
%
%
\section{Some new representations of $Q_{n}^{m}(z)$ with
\mbox{$0\leqslant m\leqslant n$}} 
\label{IV} 
\setcounter{equation}{0}
The associated Legendre function of the second kind,
$Q_{\nu}^{\mu}(z)$, may be defined \cite{Hobs31} as the following
linear combination of the Legendre functions of the first kind
$P_{\nu}^{\mu}(z)$ and $P_{\nu}^{\mu}(-z)$:
\begin{equation}
Q_{\nu}^{\mu}(z)=\frac{\pi}{2}\mathrm{e}^{\mathrm{i}\pi\mu}
\frac{\mathrm{e}^{\mp\mathrm{i}\pi\nu}P_{\nu}^{\mu}(z)
-P_{\nu}^{\mu}(-z)}{\sin[\pi(\nu+\mu)]}
\qquad (\mathrm{Im}(z)\gtrless0).
\label{4.1}
\end{equation}
In the special case of $\mu=m$, Eq.\ (\ref{4.1}) simplifies to
\begin{equation}
Q_{\nu}^{m}(z)=\frac{\pi}{2}
\frac{\mathrm{e}^{\mp\mathrm{i}\pi\nu}P_{\nu}^{m}(z)
-P_{\nu}^{m}(-z)}{\sin(\pi\nu)}
\qquad (\mathrm{Im}(z)\gtrless0),
\label{4.2}
\end{equation}
hence, after exploiting the l'Hospital rule, one obtains
\begin{equation}
Q_{n}^{m}(z)=\mp\frac{1}{2}\mathrm{i}\pi P_{n}^{m}(z)
+\frac{1}{2}\frac{\partial P_{\nu}^{m}(z)}
{\partial\nu}\bigg|_{\nu=n}
-\frac{(-)^{n}}{2}\frac{\partial P_{\nu}^{m}(-z)}
{\partial\nu}\bigg|_{\nu=n}
\qquad (\textrm{$0\leqslant m\leqslant n$,
$\mathrm{Im}(z)\gtrless0$}).
\label{4.3}
\end{equation}
Thus, we see that the problem of evaluation of $Q_{n}^{m}(z)$ with
$0\leqslant m\leqslant n$ may be reduced to that of derivation of
expressions for $[\partial P_{\nu}^{m}(\pm z)/\partial\nu]_{\nu=n}$.

Accordingly, after combining Eq.\ (\ref{4.3}) with Eqs.\ (\ref{1.3})
and (\ref{3.1}), we obtain
\begin{eqnarray}
Q_{n}^{m}(z) &=& \frac{1}{2}P_{n}^{m}(z)\ln\frac{z+1}{z-1}
\mp\frac{1}{2}[\psi(n+m+1)+\psi(n+1)]P_{n}^{m}(z)
\nonumber \\
&& \pm\,\frac{(\pm)^{n}(\mp)^{m}}{2}\frac{(n+m)!}{(n-m)!}
\left(\frac{z^{2}-1}{4}\right)^{-m/2}
\sum_{k=0}^{m-1}(\mp)^{k}\frac{(k+n-m)!(m-k-1)!}{k!(n+m-k)!}
\left(\frac{z\mp1}{2}\right)^{k}
\nonumber \\
&& \pm\,\frac{(\pm)^{n+m}}{2}\left(\frac{z^{2}-1}{4}\right)^{m/2}
\sum_{k=0}^{n-m}(\pm)^{k}
\frac{(k+n+m)!\psi(k+n+m+1)}{k!(k+m)!(n-m-k)!}
\left(\frac{z\mp1}{2}\right)^{k}
\nonumber \\
&& \mp\,\frac{(\pm)^{n}}{2}\frac{(n+m)!}{(n-m)!}
\left(\frac{z\mp1}{z\pm1}\right)^{m/2}
\sum_{k=0}^{n}(\pm)^{k}\frac{(k+n)!}{k!(k+m)!(n-k)!}
\nonumber \\
&& \quad \times[\psi(k+n+1)-\psi(k+m+1)-\psi(k+1)]
\left(\frac{z\mp1}{2}\right)^{k}
\qquad (0\leqslant m\leqslant n),
\label{4.4}
\end{eqnarray}
where the upper signs follow if $[\partial
P_{\nu}^{m}(z)/\partial\nu]_{\nu=n}$ is evaluated from Eq.\
(\ref{1.3}) and $[\partial P_{\nu}^{m}(-z)/\partial\nu]_{\nu=n}$ from
Eq.\ (\ref{3.1}), while the lower signs result if the roles of Eqs.\
(\ref{1.3}) and (\ref{3.1}) are interchanged. The same expression for
$Q_{n}^{m}(z)$ as above is obtained if Eq.\ (\ref{4.3}) is coupled
with Eqs.\ (\ref{1.4}) and (\ref{3.7}). Further, using Eqs.\
(\ref{1.3}) and (\ref{3.2}) in Eq.\ (\ref{4.3}) leads to
\begin{eqnarray}
Q_{n}^{m}(z) &=& \frac{1}{2}P_{n}^{m}(z)\ln\frac{z+1}{z-1}
\mp\frac{1}{2}[\psi(n+1)+\psi(n-m+1)]P_{n}^{m}(z)
\nonumber \\
&& \pm\,\frac{(\pm)^{n}(-)^{m}}{2}
\left(\frac{z\pm1}{z\mp1}\right)^{m/2}
\sum_{k=0}^{m-1}(\mp)^{k}\frac{(k+n)!(m-k-1)!}{k!(n-k)!}
\left(\frac{z\mp1}{2}\right)^{k}
\nonumber \\
&& \mp\,\frac{(\pm)^{n+m}}{2}\left(\frac{z^{2}-1}{4}\right)^{m/2}
\sum_{k=0}^{n-m}(\pm)^{k}\frac{(k+n+m)!}{k!(k+m)!(n-m-k)!}
\nonumber \\
&& \quad \times[\psi(k+n+m+1)-\psi(k+m+1)-\psi(k+1)]
\left(\frac{z\mp1}{2}\right)^{k}
\nonumber \\
&& \pm\,\frac{(\pm)^{n}}{2}\frac{(n+m)!}{(n-m)!}
\left(\frac{z\mp1}{z\pm1}\right)^{m/2}
\sum_{k=0}^{n}(\pm)^{k}\frac{(k+n)!\psi(k+n+1)}{k!(k+m)!(n-k)!}
\left(\frac{z\mp1}{2}\right)^{k}
\nonumber \\
&& (0\leqslant m\leqslant n).
\label{4.5}
\end{eqnarray}
Next, if $[\partial P_{\nu}^{m}(z)/\partial\nu]_{\nu=n}$ is obtained
from Eq.\ (\ref{3.1}) and $[\partial
P_{\nu}^{m}(-z)/\partial\nu]_{\nu=n}$ from Eq.\ (\ref{3.7}), or vice
versa, then Eq.\ (\ref{4.3}) yields the expressions
\begin{eqnarray}
Q_{n}^{m}(z) &=& \frac{1}{2}P_{n}^{m}(z)\ln\frac{z+1}{z-1}
\mp\psi(n+m+1)P_{n}^{m}(z)
\nonumber \\
&& \pm\,\frac{(\pm)^{n}(\mp)^{m}}{2}\frac{(n+m)!}{(n-m)!}
\left(\frac{z^{2}-1}{4}\right)^{-m/2}
\sum_{k=0}^{m-1}(\mp)^{k}\frac{(k+n-m)!(m-k-1)!}{k!(n+m-k)!}
\left(\frac{z\mp1}{2}\right)^{k}
\nonumber \\
&& \pm\,\frac{(\pm)^{n+m}}{2}\left(\frac{z^{2}-1}{4}\right)^{m/2}
\sum_{k=0}^{n-m}(\pm)^{k}\frac{(k+n+m)!}{k!(k+m)!(n-m-k)!}
\nonumber \\
&& \quad \times[2\psi(k+n+m+1)-\psi(k+m+1)]
\left(\frac{z\mp1}{2}\right)^{k}
\nonumber \\
&& \mp\,\frac{(\pm)^{n}}{2}\frac{(n+m)!}{(n-m)!}
\left(\frac{z\mp1}{z\pm1}\right)^{m/2}
\sum_{k=0}^{n}(\pm)^{k}\frac{(k+n)!}{k!(k+m)!(n-k)!}
\nonumber \\
&& \quad \times[2\psi(k+n+1)-2\psi(k+m+1)-\psi(k+1)]
\left(\frac{z\mp1}{2}\right)^{k}
\qquad (0\leqslant m\leqslant n).
\label{4.6}
\end{eqnarray}
We are not aware of any appearance of either of the formulas
(\ref{4.4})--(\ref{4.6}) in the literature. Furthermore, if Eqs.\
(\ref{3.1}) and (\ref{3.8}) are used in Eq.\ (\ref{4.3}), this
results in
\begin{eqnarray}
Q_{n}^{m}(z) &=& \frac{1}{2}P_{n}^{m}(z)\ln\frac{z+1}{z-1}
\mp\psi(n+1)P_{n}^{m}(z)
\nonumber \\
&& \pm\,\frac{(\pm)^{n}(\mp)^{m}}{2}\frac{(n+m)!}{(n-m)!}
\left(\frac{z^{2}-1}{4}\right)^{-m/2}
\sum_{k=0}^{m-1}(\mp)^{k}\frac{(k+n-m)!(m-k-1)!}{k!(n+m-k)!}
\left(\frac{z\mp1}{2}\right)^{k}
\nonumber \\
&& \pm\,\frac{(\pm)^{n+m}}{2}\left(\frac{z^{2}-1}{4}\right)^{m/2}
\sum_{k=0}^{n-m}(\pm)^{k}\frac{(k+n+m)!\psi(k+m+1)}{k!(k+m)!(n-m-k)!}
\left(\frac{z\mp1}{2}\right)^{k}
\nonumber \\
&& \pm\,\frac{(\pm)^{n}}{2}\frac{(n+m)!}{(n-m)!}
\left(\frac{z\mp1}{z\pm1}\right)^{m/2}
\sum_{k=0}^{n}(\pm)^{k}\frac{(k+n)!\psi(k+1)}{k!(k+m)!(n-k)!}
\left(\frac{z\mp1}{2}\right)^{k}
\nonumber \\
&& (0\leqslant m\leqslant n),
\label{4.7}
\end{eqnarray}
which is the same what follows if Eqs.\ (\ref{1.3}) and (\ref{1.4})
are plugged into Eq.\ (\ref{4.3}) (cf.\ Ref.\ \cite{Szmy09b}).
Finally, insertion of Eqs.\ (\ref{3.2}) and (\ref{3.7}) into Eq.\
(\ref{4.3}) leads to
\begin{eqnarray}
Q_{n}^{m}(z) &=& \frac{1}{2}P_{n}^{m}(z)\ln\frac{z+1}{z-1}
\mp\frac{1}{2}[\psi(n+m+1)+\psi(n-m+1)]P_{n}^{m}(z)
\nonumber \\
&& \pm\,\frac{(\pm)^{n}(-)^{m}}{2}
\left(\frac{z\pm1}{z\mp1}\right)^{m/2}
\sum_{k=0}^{m-1}(\mp)^{k}\frac{(k+n)!(m-k-1)!}{k!(n-k)!}
\left(\frac{z\mp1}{2}\right)^{k}
\nonumber \\
&& \pm\,\frac{(\pm)^{n+m}}{2}\left(\frac{z^{2}-1}{4}\right)^{m/2}
\sum_{k=0}^{n-m}(\pm)^{k}\frac{(k+n+m)!\psi(k+1)}{k!(k+m)!(n-m-k)!}
\left(\frac{z\mp1}{2}\right)^{k}
\nonumber \\
&& \pm\,\frac{(\pm)^{n}}{2}\frac{(n+m)!}{(n-m)!}
\left(\frac{z\mp1}{z\pm1}\right)^{m/2}
\sum_{k=0}^{n}(\pm)^{k}\frac{(k+n)!\psi(k+m+1)}{k!(k+m)!(n-k)!}
\left(\frac{z\mp1}{2}\right)^{k}
\nonumber \\
&& (0\leqslant m\leqslant n),
\label{4.8}
\end{eqnarray}
which, in turn, is the same what is obtained if Eqs.\ (\ref{1.3}) and
(\ref{1.5}) are coupled with Eq.\ (\ref{4.3}) (cf.\ again Ref.\
\cite{Szmy09b}; for alternative derivations of the above result see
Ref.\ \cite[pages 81, 82 and 85]{Robi58} and Ref.\ \cite{Szmy09a}).
Other expressions for $Q_{n}^{m}(z)$ may be obtained by combining
Eqs.\ (\ref{4.4})--(\ref{4.8}), with the possible help of the
identities (\ref{3.4}) and (\ref{3.6}).

Once the function $Q_{n}^{m}(z)$ is known, one may find the function
$Q_{n}^{-m}(z)$ from the well-known relationship
\begin{equation}
Q_{n}^{-m}(z)=\frac{(n-m)!}{(n+m)!}Q_{n}^{m}(z)
\qquad (0\leqslant m\leqslant n).
\label{4.9}
\end{equation}

Using the formula
\begin{equation}
Q_{n}^{\pm m}(x)=\frac{(-)^{m}}{2}
\left[\mathrm{e}^{\mp\mathrm{i}m/2}Q_{n}^{\pm m}(x+\mathrm{i}0)
+\mathrm{e}^{\pm\mathrm{i}m/2}Q_{n}^{\pm m}(x-\mathrm{i}0)\right],
\label{4.10}
\end{equation}
which follows from the Hobson's \cite{Hobs31} definition of the
associated Legendre function of the second kind on the cut
$x\in[-1,1]$, and the identities in Eq.\ (\ref{3.13}), from Eqs.\
(\ref{4.4})--(\ref{4.9}) one may derive counterpart representations
of $Q_{n}^{\pm m}(x)$ with $0\leqslant m\leqslant n$. Since the
procedure does not offer any difficulty, we do not present the
resulting expressions here.
%
%

%

\begin{thebibliography}{99}
\bibitem{Bryc09}
   Yu.\ A.\ Brychkov,
   On the derivatives of the Legendre functions $P_{\nu}^{\mu}(z)$
   and $Q_{\nu}^{\mu}$ with respect to $\mu$ and $\nu$,
   Integral Transforms Spec.\ Funct.\ doi:10.1080/10652460903069660
\bibitem{Cohl09}
   H.\ S.\ Cohl,
   Derivatives with respect to the degree and order of associated
   Legendre functions for $|z|>1$ using modified Bessel functions,
   Integral Transforms Spec.\ Funct.\ in press
\bibitem{Szmy09a}
   R.\ Szmytkowski,
   On the derivative of the associated Legendre function of the 
   first kind of integer degree with respect to its order (with
   applications to the construction of the associated Legendre
   function of the second kind of integer degree and order),
   J.\ Math.\ Chem.\ 46 (2009) 231
\bibitem{Szmy09b}
   R.\ Szmytkowski,
   On the derivative of the associated Legendre function of the first
   kind of integer order with respect to its degree (with
   applications to the construction of the associated Legendre
   function of the second kind of integer degree and order),
   preprint arXiv:0907.3217
\bibitem{Hobs31}
   E.\ W.\ Hobson,
   The Theory of Spherical and Ellipsoidal Harmonics,
   Cambridge University Press, Cambridge, 1931
   [reprinted: (Chelsea, New York, 1955)]
\bibitem{Robi57}
   L.\ Robin,
   Fonctions Sph{\'e}riques de Legendre et Fonctions
   Sph{\'e}ro{\"{\i}}dales, vol.\ 1,
   Gauthier-Villars, Paris, 1957
\bibitem{Robi58}
   L.\ Robin,
   Fonctions Sph{\'e}riques de Legendre et Fonctions
   Sph{\'e}ro{\"{\i}}dales, vol.\ 2,
   Gauthier-Villars, Paris, 1958
\bibitem{Robi59}
   L.\ Robin,
   Fonctions Sph{\'e}riques de Legendre et Fonctions
   Sph{\'e}ro{\"{\i}}dales, vol.\ 3,
   Gauthier-Villars, Paris, 1959
\bibitem{Mors53}
   P.\ M.\ Morse, H.\ Feshbach,
   Methods of Theoretical Physics,
   McGraw-Hill, New York, 1953
\bibitem{Moon61}
   P.\ Moon, D.\ E.\ Spencer,
   Field Theory for Engineers,
   Van Nostrand, Princeton, 1961
\bibitem{Barn07}
   E.\ W.\ Barnes,
   On generalized Legendre functions,
   Quart.\ J.\ Pure Appl.\ Math.\ 39 (1907) 97
\bibitem{Szmy06}
   R.\ Szmytkowski,
   On the derivative of the Legendre function of the first kind 
   with respect to its degree,
   J.\ Phys.\ A 39 (2006) 15147
   [erratum: 40 (2007) 7819]
\bibitem{Szmy07}
   R.\ Szmytkowski, 
   Addendum to `On the derivative of the Legendre function of the
   first kind with respect to its degree,'
   J.\ Phys.\ A 40 (2007) 14887.
   We note parenthetically that the simplest way to arrive at Eq.\ 
   (5) in that paper is to differentiate both sides of the identity
   $\mathrm{d}^{k}z^{\alpha}/\mathrm{d}z^{k}
   =[\Gamma(\alpha+1)/\Gamma(\alpha-k+1)]z^{\alpha-k}$ with respect
   to $\alpha$.
\bibitem{Sche41}
   S.\ A.\ Schelkunoff,
   Theory of antennas of arbitrary size and shape,
   Proc.\ IRE 29 (1941) 493
   [reprinted: Proc.\ IEEE 72 (1984) 1165].
   Notice that Schelkunoff used a definition of the digamma function
   different from that adopted in Ref.\ \cite{Szmy06} and in the 
   present paper.
\end{thebibliography}
\end{document}